\documentstyle[12pt]{article}

\font\twlgot =eufm10 scaled \magstep1
\font\egtgot =eufm8
\font\sevgot =eufm7
\font\twlmsb =msbm10 scaled \magstep1
\font\egtmsb =msbm8
\font\sevmsb =msbm7

\newfam\gotfam
\def\pgot{\fam\gotfam\twlgot}
\textfont\gotfam\twlgot
\scriptfont\gotfam\egtgot
\scriptscriptfont\gotfam\sevgot
\def\got{\protect\pgot}
\newfam\msbfam
\textfont\msbfam\twlmsb
\scriptfont\msbfam\egtmsb
\scriptscriptfont\msbfam\sevmsb
\def\Bbb{\protect\pBbb}
\def\pBbb{\relax\ifmmode\expandafter\Bb\else\typeout{You cann't use
Bbb in text mode}\fi}
\def\Bb #1{{\fam\msbfam\relax#1}}

\newcommand{\gO}{{\got O}}
\newcommand{\gQ}{{\got T}}

\newcommand{\gA}{{\got A}}

\newcommand{\gd}{{\got d}}

\newcommand{\gS}{{\got S}}

\def\thebibliography#1{\section*{References}\list
  {[\arabic{enumi}]}{\settowidth\labelwidth{#1}\leftmargin\labelwidth
    \advance\leftmargin\labelsep
    \usecounter{enumi}}
    \def\newblock{\hskip .11em plus .33em minus .07em}
    \sloppy\clubpenalty4000\widowpenalty4000
    \sfcode`\.=1000\relax}

\def\op#1{\mathop{\fam0 #1}\limits}

\newcommand{\nm}[1]{|{#1}|}
\newcommand{\beq}{\begin{equation}}
\newcommand{\eeq}{\end{equation}}
\newcommand{\ben}{\begin{eqnarray}}
\newcommand{\een}{\end{eqnarray}}
\newcommand{\be}{\begin{eqnarray*}}
\newcommand{\ee}{\end{eqnarray*}}
\newcommand{\bea}{\begin{eqalph}}
\newcommand{\eea}{\end{eqalph}}
\newcommand{\cA}{{\cal A}}

\newcommand{\cP}{{\cal P}}

\newcommand{\cL}{{\cal L}}

\newcommand{\cE}{{\cal E}}

\newcommand{\cS}{{\cal S}}

\newcommand{\cO}{{\cal O}}

\newcommand{\cK}{{\cal K}}
\newcommand{\bL}{{\bf L}}

\newcommand{\bE}{{\bf E}}

\newcommand{\vr}{\varrho}

\newcommand{\dl}{\delta}
\newcommand{\la}{\lambda}
\newcommand{\La}{\Lambda}
\newcommand{\f}{\phi}
\newcommand{\om}{\omega}
\newcommand{\Om}{\Omega}
\newcommand{\m}{\mu}

\newcommand{\G}{\Gamma}
\newcommand{\th}{\theta}

\newcommand{\vt}{\vartheta}
\newcommand{\vf}{\varphi}

\newcommand{\di}{{\rm dim\,}}

\newcommand{\si}{\sigma}
\newcommand{\Si}{\Sigma}
\newcommand{\w}{\wedge}

\newcommand{\ol}{\overline}
\newcommand{\dr}{\partial}
\newcommand{\ar}{\op\longrightarrow}

\let\ssection=\section
\renewcommand{\section}{\setcounter{equation}{0}\ssection}

\newcounter{eqalph}
\newcounter{equationa}
\newcounter{remark}
\newcounter{example}
\newcounter{theorem}
\newcounter{proposition}
\newcounter{lemma}
\newcounter{corollary}
\newcounter{definition}
\def\theremark{\arabic{remark}}

\def\thedefinition{\arabic{definition}}

\newenvironment{proof}{\noindent
{\it Proof.}}{$\Box$ \medskip}
\newenvironment{rem}{\refstepcounter{remark}\medskip\noindent{\it
Remark \theremark.}}{\medskip}
\newenvironment{theo}{\refstepcounter{definition}
\bigskip\noindent{\bf Theorem \thedefinition.} \it}{\medskip}
\newenvironment{prop}{\refstepcounter{definition}
\bigskip\noindent{\bf Proposition \thedefinition.}\it}{\medskip}
\newenvironment{lem}{\refstepcounter{definition}
\bigskip\noindent{\bf Lemma \thedefinition.}\it}{\medskip}

\newenvironment{eqalph}{\stepcounter{equation}
\setcounter{equationa}{\value{equation}} \setcounter{equation}{0}

\begin{eqnarray}}{\end{eqnarray}\setcounter{equation}{\value{equationa}}}

\newcommand{\mar}[1]{}

\begin{document}
\hbox{}

{\parindent=0pt

{\large \bf The variational bicomplex on graded manifolds and
its cohomology}
\bigskip

{\bf G. Sardanashvily}
\bigskip

\begin{small}

Department of Theoretical Physics, 
Moscow State University, 117234 Moscow, Russia

\bigskip

{\bf Abstract:} Lagrangian formalism on graded manifolds is phrased in
terms of the Grassmann-graded variational bicomplex, generalizing the
familiar variational bicomplex for even Lagrangian systems on fiber
bundles.

\end{small}

  }

\bigskip
\bigskip

Lagrangian systems of odd and affine even fields on a smooth manifold $X$
($\di X=n$) can be described in algebraic terms of the
Grassmann-graded variational bicomplex
\cite{barn,jmp05,cmp04}, generalizing the variational
bicomplex for even Lagrangian systems on fiber bundles
\cite{ander,jmp,tak2}. Here, this bicomplex is stated n a
general setting, when a fiber bundle $Y\to X$ of even fields need not be
affine. For this purpose, we consider graded manifolds whose body is a
fiber bundle
$Y\to X$ and its jet manifolds $J^rY$, but not
$X$. We show that the relevant cohomology of the Grassmann-graded
variational bicomplex on these graded manifolds reduces to that of the
variational bicomplex on a fiber bundle $Y\to X$. 

\begin{rem} 
Smooth
manifolds throughout are assumed to be real,
finite-dimensional, Hausdorff, second-countable (consequently,
paracompact) and connected.
By a Grassmann algebra
over a ring $\cK$ is meant a $\Bbb Z_2$-graded exterior algebra of some
$\cK$-module. 
We restrict our consideration to graded manifolds $(Z,\gA)$
with structure sheaves
$\gA$ of Grassmann algebras of finite rank
\cite{bart,book05}. The symbols $|.|$ and $[.]$ stand for the form degree
and Grassmann parity, respectively. We denote by $\La$, $\Si$, $\Xi$,
$\Om$ the symmetric multi-indices, e.g.,  $\La=(\la_1...\la_k)$, 
$\la+\La=(\la\la_1...\la_k)$. Summation over
a multi-index
$\La=(\la_1...\la_k)$ throughout means separate summation over each its
index $\la_i$.
\end{rem}

Let $J^rY$, $r\in\Bbb N$, be finite order jet
manifolds of sections of
$Y\to X$, where $r=0$ conventionally stands for $Y$. They
make up the inverse system
\mar{5.10}\beq
X\op\longleftarrow^\pi Y\op\longleftarrow^{\pi^1_0} J^1Y
\longleftarrow \cdots J^{r-1}Y \op\longleftarrow^{\pi^r_{r-1}}
J^rY\longleftarrow\cdots, \label{5.10}
\eeq
where  $\pi^r_{r-1}$ are affine bundles and, hence, open maps. Its
projective limit
$(J^\infty Y,
\pi^\infty_r:J^\infty Y\to J^rY)$ is a paracompact
Fr\'echet manifold, called the infinite order jet
manifold.  Moreover, $Y$ is the
strong deformation retract of $J^\infty Y$. A bundle atlas
$\{(U;x^\la,y^i)\}$ of
$Y\to X$ yields the coordinate atlas
\mar{jet1}\beq
\{((\pi^\infty_0)^{-1}(U); x^\la, y^i_\La)\}, \qquad
{y'}^i_{\la+\La}=\frac{\dr x^\m}{\dr x'^\la}d_\m y'^i_\La, \qquad
0\leq|\La|, \label{jet1}
\eeq
of $J^\infty Y$, where 
\mar{5.177}\beq
d_\la = \dr_\la + \op\sum_{0\leq|\La|} y^i_{\la+\La}\dr_i^\La,
\qquad d_\La=d_{\la_1}\circ\cdots\circ d_{\la_k},  \label{5.177}
\eeq
are total derivatives. Let us fix an atlas of
$Y$  containing a finite number of charts \cite{greub}.

The inverse system (\ref{5.10}) yields the direct system
\mar{5.7}\beq
\cO^*X\op\longrightarrow^{\pi^*} \cO^*Y
\op\longrightarrow^{\pi^1_0{}^*} \cO_1^*Y \ar\cdots \cO^*_{r-1}Y
\op\longrightarrow^{\pi^r_{r-1}{}^*}
 \cO_r^*Y \longrightarrow\cdots,  \label{5.7}
\eeq
of graded differential algebras (henceforth GDAs)  $\cO_r^*Y$ of
exterior forms on jet manifolds $J^rY$ with respect to the
pull-back monomorphisms $\pi^r_{r-1}{}^*$. Its direct limit
is the GDA $\cO_\infty^*Y$
of all exterior forms on finite order jet manifolds
modulo the pull-back identification.
 One can think of elements of
$\cO_\infty^*Y$ as being exterior forms on the infinite order jet manifold
$J^\infty Y$ as follows. Let $\gO^*_r$ be the sheaf
of germs of exterior forms on  $J^rY$ and 
$\ol\gO^*_r$ the canonical presheaf of local sections of $\gO^*_r$,
seen as a particular topological bundle over $Y$ (we follow the
terminology of
\cite{hir}). Since
$\pi^r_{r-1}$ are open maps, there is the direct  system
of presheaves
\be
\ol\gO^*_X\op\longrightarrow^{\pi^*} \ol\gO^*_0
\op\longrightarrow^{\pi^1_0{}^*} \ol\gO_1^* \cdots
\op\longrightarrow^{\pi^r_{r-1}{}^*}
 \ol\gO_r^* \longrightarrow\cdots.
\ee
Its direct limit $\ol\gO^*_\infty$
is a presheaf of GDAs on 
$J^\infty Y$. Let
$\gQ^*_\infty$ be the sheaf of GDAs of germs of $\ol\gO^*_\infty$ on
$J^\infty Y$. The structure module
$\G(\gQ^*_\infty)$ of global sections of $\gQ^*_\infty$ is a GDA such
that, given an element
$\f\in \G(\gQ^*_\infty)$ and a point $z\in J^\infty Y$, there
exist an open
neighbourhood $U$ of $z$ and an
exterior form
$\f^{(k)}$ on some finite order jet manifold $J^kY$ so that
$\f|_U= \pi^{\infty*}_k\f^{(k)}|_U$. Therefore, there is the GDA
monomorphism 
\mar{c0}\beq
\cO^*_\infty Y
\to\G(\gQ^*_\infty). \label{c0}
\eeq 
It should be emphasized that the paracompact space $J^\infty Y$ admits
a partition of
unity by elements of the ring $\G(\gQ^0_\infty)$.

Due to the monomorphism (\ref{c0}), one can restrict 
$\cO^*_\infty Y$ to the coordinate chart (\ref{jet1}) where horizontal
forms
$\{dx^\la\}$ and contact one-forms $\{\th^i_\La=dy^i_\La
-y^i_{\la+\La}dx^\la\}$ make up a local basis for the
$\cO^0_\infty Y$-algebra $\cO^*_\infty Y$. Though $J^\infty Y$
is not a smooth manifold, elements of $\cO^*_\infty Y$ are
exterior forms on finite order jet manifolds and, therefore, their
coordinate transformations are smooth.
There is the canonical decomposition $\cO^*_\infty
Y=\oplus\cO^{k,m}_\infty Y$ of $\cO^*_\infty Y$ into $\cO^0_\infty
Y$-modules $\cO^{k,m}_\infty Y$ of $k$-contact and $m$-horizontal
forms together with the corresponding projectors 
\be
h_k:\cO^*_\infty
Y\to \cO^{k,*}_\infty Y, \qquad h^m:\cO^*_\infty Y\to
\cO^{*,m}_\infty Y.
\ee
Accordingly, the exterior differential on
$\cO_\infty^* Y$ is split into the sum $d=d_H+d_V$ of the
nilpotent total and vertical differentials, where
\be
d_H\circ h_k=h_k\circ d\circ h_k, \qquad d_H\circ
h_0=h_0\circ d, \qquad d_H(\f)= dx^\la\w d_\la(\f).
\ee
One also introduces the $\Bbb R$-module projector
\be
\vr: \cO^{k,n}_\infty Y\to \bE_k\subset \cO^{k,n}_\infty Y, \qquad k=1,
\ldots,  
\ee
such that $\vr\circ d_H=0$ and the nilpotent
variational operator $\dl=\vr\circ d$ on $\cO^{*,n}_\infty Y$. 
Then the GDA $\cO^*_\infty
Y$ is split into the above mentioned variational bicomplex. This
contains the variational subcomplex
\mar{b317}\beq
0\to\Bbb R\to \cO^0_\infty Y \ar^{d_H}\cO^{0,1}_\infty Y\cdots
\op\ar^{d_H} \cO^{0,n}_\infty Y \op\ar^\dl \bE_1 \op\ar^\dl \bE_2
\ar \cdots, \label{b317}
\eeq
whose elements $L\in \cO^{0,n}_\infty Y$ and $\dl L\in \bE_1$ 
are respectively finite order Lagrangians and their Euler--Lagrange
operators on a fiber bundle $Y\to X$.

Turn now to Lagrangian systems both of even and odd fields.
Though there are different approaches to treat odd fields on a smooth
manifold $X$,  the following variant of the Serre--Swan theorem
motivates us to describe them in
terms of graded manifolds whose body is $X$. 

\begin{theo} \label{v0} \mar{v0}
Let $Z$ be a smooth manifold. A Grassmann algebra $\cA$ over the
ring $C^\infty(Z)$ of smooth real functions on $Z$ is isomorphic to the
Grassmann algebra of graded functions on a graded manifold with a body
$Z$ iff it is the exterior algebra 
of some projective $C^\infty(Z)$-module of finite rank.
\end{theo}

\begin{proof} The proof follows at once from Batchelor's theorem
\cite{bart} and the Serre-Swan theorem generalized to an arbitrary smooth
manifold \cite{book05,ren}.  
By virtue of the first one,
any graded manifold
$(Z,\gA)$ with a body
$Z$ is isomorphic to the one
$(Z,\gA_Q)$ with the structure sheaf $\gA_Q$ of germs of sections of
the exterior bundle product
\mar{g80}\beq
\w Q^*=\Bbb R\op\oplus_Z Q^*\op\oplus_Z\op\w^2
Q^*\op\oplus_Z\cdots, \label{g80}
\eeq
where $Q^*$ is the dual of some vector bundle $Q\to Z$.
We agree to call $(Z,\gA_Q)$ the simple graded manifold modelled over
the structure vector bundle $Q\to X$. Its structure ring $\cA_Q$ of graded
functions (sections of
$\gA_Q$) consists of sections of the exterior bundle
(\ref{g80}), and it is the $\Bbb Z_2$-graded exterior algebra of the
$C^\infty(Z)$-module of sections of
$Q^*\to Z$. The Serre--Swan theorem states that a
$C^\infty(Z)$-module is isomorphic to the module of sections of a smooth
vector bundle over $Z$ iff it is a projective module of finite rank.
\end{proof}

In field models, Batchelor's isomorphism is
usually fixed from the beginning. Therefore, we restrict our consideration
to simple graded manifolds $(Z,\gA_Q)$. One associates to 
$(Z,\gA_Q)$ the following bigraded differential algebra (henceforth BGDA)
$\cS^*[Q;Z]$
\cite{bart,book05}.  Let us consider the sheaf $\gd\gA_Q$ of
graded derivations of $\gA_Q$.
One can show that its sections over an open subset $U\subset Z$ exhaust
all 
$\Bbb Z_2$-graded derivations of the $\Bbb Z_2$-graded $\Bbb R$-ring
$\cA_U$ of graded functions on $U$
\cite{bart}. Global sections of $\gd\cA_Q$ make up the real Lie
superalgebra of $\Bbb Z_2$-graded derivations
of the $\Bbb R$-ring $\cA_Q$.
Then one can construct the Chevalley--Eilenberg complex of 
$\gd\cA_Q$ with coefficients in $\cA_Q$ 
\cite{fuks}. Its subcomplex $\cS^*[Q;Z]$ of $\cA_Q$-linear morphism is
the $\Bbb Z_2$-graded Chevalley--Eilenberg differential
calculus 
\mar{v1}\beq
0\to \Bbb R\to \cA_Q \ar^d \cS^1[Q;Z]\ar^d\cdots \cS^k[Q;Z]\ar^d\cdots
\label{v1}
\eeq
over a 
$\Bbb Z_2$-graded commutative $\Bbb R$-ring
$\cA_Q$ \cite{book05}. The Chevalley--Eilenberg coboundary
operator $d$ and the graded exterior product $\w$ 
make $\cS^*[Q;Z]$ into a BGDA whose elements obey the relations 
\mar{v21,2}\ben
&&\f\w\f' =(-1)^{|\f||\f'| +[\f][\f']}\f'\w \f, \label{v21}\\
&& d(\f\w\f')= d\f\w\f' +(-1)^{|\f|}\f\w d\f'.  \label{v22}
\een  
Given the GDA $\cO^*Z$ of
exterior forms on
$Z$, there are the canonical monomorphism $\cO^*Z\to \cS^*[Q;Z]$ and
body epimorphism $\cS^*[Q;Z]\to \cO^*Z$. 

\begin{lem} \label{v62} \mar{v62} 
The BGDA $\cS^*[Q;Z]$ is a minimal differential calculus over $\cA_Q$,
i.e., it is generated by elements
$df$, $f\in
\cA_Q$.
\end{lem}

\begin{proof}
One can show that
elements of
$\gd\cA_Q$ are represented by sections of some vector bundle over
$Z$, i.e., 
$\gd\cA_Q$ is a projective $C^\infty(Z)$- and $\cA_Q$-module of
finite rank, and so is its $\cA_Q$-dual $\cS^1[Q;Z]$
\cite{cmp04,book05}. Hence, 
$\gd\cA_Q$ is the $\cA_Q$-dual of $\cS^1[Q;Z]$ and, consequently, 
$\cS^1[Q;Z]$ is generated by elements
$df$, $f\in
\cA_Q$  \cite{book05}. 
\end{proof}

This fact is essential for our consideration because of the following
\cite{book05}.

\begin{lem} \label{v30} \mar{v30}
Given a ring $R$, let $\cK$, $\cK'$ be $R$-rings and  $\cA$, $\cA'$ the
Grassmann algebras over  $\cK$ and $\cK'$, respectively.
Then any homomorphism $\rho:
\cA\to \cA'$ yields the homomorphism of the minimal
Chevalley--Eilenberg differential calculus over a $\Bbb Z_2$-graded
$R$-ring $\cA$ to that over
$\cA'$ given by the map $da \mapsto d(\rho(a))$, $a\in\cA$. This map
provides a monomorphism if $\rho$ is a monomorphism of $R$-algebras
\end{lem}

One can think of elements of the BGDA $\cS^*[Q;Z]$ as being
graded exterior forms on $Z$ as follows. 
Given an open subset $U\subset Z$, let $\cA_U$ be the Grassmann algebra of
sections of the sheaf $\gA_Q$ over $U$, and let 
$\cS^*[Q;U]$ be the corresponding
Chevalley--Eilenberg differential calculus over $\cA_U$. Given an open set
$U'\subset U$, the restriction morphisms $\cA_U\to\cA_{U'}$ 
yield the restriction
morphism  of the BGDAs
$\cS^*[Q;U]\to \cS^*[Q;U']$. Thus,  we obtain
the presheaf 
$\{U,\cS^*[Q;U]\}$ of BGDAs on a manifold $Z$ and the sheaf $\gS^*[Q;Z]$
of BGDAs of germs of this presheaf. Since $\{U,\cA_U\}$ is the canonical
presheaf of the sheaf
$\gA_Q$, the
canonical presheaf of 
$\gS^*[Q;Z]$ is $\{U,\cS^*[Q;U]\}$. In particular, 
$\cS^*[Q;Z]$ is the BGDA of global sections of the sheaf $\gS^*[Q;Z]$,
and  there is the restriction morphism $\cS^*[Q;Z]\to \cS^*[Q;U]$ for
any open $U\subset Z$. Due to this morphism, elements of $\cS^*[Q;Z]$
can be written in the following local form.

Given bundle
coordinates
$(z^A,q^a)$ on
$Q$ and the
corresponding fiber basis $\{c^a\}$ for $Q^*\to X$, the tuple
$(z^A, c^a)$ is called a local basis for the graded manifold
$(Z,\gA_Q)$ \cite{bart}. With respect to this basis, graded functions read
\mar{v23}\beq
f=\op\sum_{k=0} \frac1{k!}f_{a_1\ldots a_k}c^{a_1}\cdots
c^{a_k}, \label{v23}
\eeq
where $f_{a_1\cdots a_k}$ are smooth real functions on $Z$,
and we omit the symbol of the exterior product of elements $c^a$.
Due to the canonical splitting $VQ= Q\times Q$, 
the fiber basis $\{\dr_a\}$ for vertical tangent bundle $VQ\to Q$ of $Q\to
Z$ is the dual of $\{c^a\}$. Then 
graded derivations take the local form $u=
u^A\dr_A + u^a\dr_a$, where $u^A, u^a$ are local graded functions.
They act on graded functions (\ref{v23}) by the rule
\mar{cmp50'}\beq
u(f_{a\ldots b}c^a\cdots c^b)=u^A\dr_A(f_{a\ldots b})c^a\cdots
c^b +u^d f_{a\ldots b}\dr_d\rfloor (c^a\cdots c^b). \label{cmp50'}
\eeq
Relative to the dual local bases
$\{dz^A\}$ for $T^*Z$ and $\{dc^b\}$ for $Q^*$, graded one-forms
read $\f=\f_A dz^A + \f_adc^a$. The duality morphism
is given by the interior product
\be
u\rfloor \f=u^A\f_A + (-1)^{[\f_a]}u^a\f_a, \qquad u\in \gd\cA_Q, \qquad
\f\in \cS^1[Q;Z]. 
\ee
The Chevalley--Eilenberg coboundary operator $d$, called the
graded exterior differential, reads 
\be
d\f=dz^A\w \dr_A\f + dc^a\w \dr_a\f,
\ee
where the derivations $\dr_A$ and $\dr_a$ act on coefficients of graded
exterior forms by the formula (\ref{cmp50'}), and they are graded
commutative with the graded exterior forms $dz^A$ and $dc^a$. 

Involving even fields which need not be affine, we come to graded
manifolds whose body is a fiber bundle $Y\to X$. We define 
jets of odd fields as simple graded manifolds
modelled over jet bundles over $X$
\cite{jmp05,cmp04}. This definition 
differs from that of jets of a graded commutative ring \cite{book05} and
jets of a graded fiber bundle
\cite{hern}, but reproduces the heuristic notion of jets of odd ghosts
 in the Lagrangian BRST theory
\cite{barn,bran01}.

Given a vector bundle $F\to X$, let us consider the simple graded manifold
$(J^rY,\gA_{F_r})$ whose body is  
$J^rY$ and structure vector bundle is the pull-back
\be
F_r=J^rY\op\times_XJ^rF
\ee
onto $J^rY$ of the jet bundle $J^rF\to X$. Given the simple
graded manifold 
$(J^{r+1}Y,\gA_{F_{r+1}})$, there is an epimorphism of graded
manifolds 
\be
(J^{r+1}Y,\gA_{F_{r+1}}) \to (J^rY,\gA_{F_r}), 
\ee
seen as local-ringed spaces. It consists of the surjection
$\pi^{r+1}_r$ and the sheaf monomorphism 
$\pi_r^{r+1*}\gA_{F_r}\to \gA_{F_{r+1}}$, 
where $\pi_r^{r+1*}\gA_{F_r}$ is the pull-back onto $J^{r+1}Y$ of the
topological fiber bundle $\gA_{F_r}\to J^rY$. This sheaf monomorphism
induces the monomorphism of the canonical presheaves 
\mar{v33}\beq
\ol \gA_{F_r}\to
\ol \gA_{F_{r+1}}, \label{v33}
\eeq
which associates to each open subset $U\subset J^{r+1}Y$ the ring of
sections of $\gA_{F_r}$ over $\pi^{r+1}_r(U)$. Accordingly, there
is the monomorphsism of
$\Bbb Z_2$-graded rings $\cA_{F_r} \to \cA_{F_{r+1}}$. 
By virtue of Lemmas \ref{v62} and \ref{v30}, this monomorphism
yields the monomorphism of BGDAs 
\mar{v4}\beq
\cS^*[F_r;J^rY]\to \cS^*[F_{r+1};J^{r+1}Y]. \label{v4}
\eeq
As a consequence, we have 
the direct system of BGDAs
\mar{v6}\beq
\cS^*[Y\op\times_X F;Y]\ar \cS^*[F_1;J^1Y]\ar\cdots 
\cS^*[F_r;J^rY]\ar\cdots, \label{v6} 
\eeq
whose direct limit $\cS^*_\infty[F;Y]$  is a BGDA of all graded
differential forms $\f\in \cS^*[F_r;J^rY]$ on jet manifolds $J^rY$ modulo
monomorphisms (\ref{v4}). Its elements obey the relations
(\ref{v21}) -- (\ref{v22}).

The monomorphisms $\cO^*_rY\to \cS^*[F_r;J^rY]$ provide
a monomorphism of the direct system (\ref{5.7}) to the direct system
(\ref{v6}) and, consequently, the monomorphism 
\mar{v7}\beq
\cO^*_\infty Y\to \cS^*_\infty[F;Y] \label{v7}
\eeq
of their direct limits. In particular, $\cS^*_\infty[F;Y]$ is an
$\cO^0_\infty Y$-algebra. Accordingly, the body epimorphisms
$\cS^*[F_r;J^rY]\to \cO^*_rY$ yield the epimorphism of $\cO^0_\infty
Y$-modules
\mar{v7'}\beq
\cS^*_\infty[F;Y]\to \cO^*_\infty Y.  \label{v7'}
\eeq

If
$Y\to X$ is an affine bundle, we recover the BGDA introduced in
\cite{jmp05,cmp04} by restricting the ring $\cO^0_\infty Y$ to its
subring $\cP^0_\infty Y$ of polynomial functions, but now one should
regard elements of $\cS^*_\infty[F;Y]$ as graded exterior forms on
the infinite order jet manifold $J^\infty Y$, but not $X$.

Indeed, let $\gS^*[F_r;J^rY]$ be the sheaf of BGDAs on $J^rY$ and
$\ol\gS^*[F_r;J^rY]$ its canonical presheaf whose elements are the
Chevalley--Eilenberg differential calculus over elements of the presheaf
$\ol\gA_{F_r}$. Then the presheaf monomorphisms (\ref{v33}) yield the
direct system of presheaves
\mar{v15}\beq
\ol\gS^*[Y\times F;Y]\ar \ol\gS^*[F_1;J^1Y] \ar\cdots
\ol\gS^*[F_r;J^rY]  \ar\cdots, \label{v15}
\eeq 
whose direct limit $\ol\gS_\infty^*[F;Y]$ is a presheaf of BGDAs on the
infinite order jet manifold $J^\infty Y$. 
Let
$\gQ^*_\infty[F;Y]$ be the sheaf of BGDAs of germs of 
the presheaf $\ol\gS_\infty^*[F;Y]$.  The structure module
$\G(\gQ^*_\infty[F;Y])$ of
sections of $\gQ^*_\infty[F;Y]$ is a BGDA such that, given an
element
$\f\in \G(\gQ^*_\infty[F;Y])$ and a point $z\in J^\infty Y$, there
exist an open
neighbourhood $U$ of $z$ and a graded
exterior form
$\f^{(k)}$ on some finite order jet manifold $J^kY$ so that
$\f|_U= \pi^{\infty*}_k\f^{(k)}|_U$. In
particular, there is the  monomorphism $\cS^*_\infty[F;Y]
\to\G(\gQ^*_\infty[F;Y])$. 

Due to this monomorphism, one can restrict 
$\cS^*_\infty[F;Y]$ to the coordinate chart (\ref{jet1}) and say 
that $\cS^*_\infty[F;Y]$ as an
$\cO^0_\infty Y$-algebra is locally generated by  the elements
\be
(1, c^a_\La, 
dx^\la,\th^a_\La=dc^a_\La-c^a_{\la+\La}dx^\la,\th^i_\La=
dy^i_\La-y^i_{\la+\La}dx^\la), \qquad 0\leq |\La|,
\ee
where $c^a_\La$, $\th^a_\La$ are odd and $dx^\la$,
$\th^i_\La$ are even. We agree to call $(y^i,c^a)$ the local basis for  
$\cS^*_\infty[F;Y]$. Let the collective symbol $s^A$ stand for
its elements. Accordingly, the notation $s^A_\La$ and
$\th^A_\La=ds^A_\La- s^A_{\la+\La}dx^\la$ is introduced. For the sake of
simplicity, we further denote $[A]=[s^A]$.

Similarly to $\cO^*_\infty Y$, the BGDA $\cS^*_\infty[F;Y]$ is
decomposed into $\cS^0_\infty[F;Y]$-modules
$\cS^{k,r}_\infty[F;Y]$ of $k$-contact and $r$-horizontal graded
forms.
Accordingly, the graded exterior differential $d$ on
$\cS^*_\infty[F;Y]$ falls into the sum $d=d_H+d_V$ of the total
and vertical differentials, where
\be
d_H(\f)=dx^\la\w d_\la(\f), \qquad d_\la = \dr_\la + \op\sum_{0\leq|\La|}
s^A_{\la+\La}\dr_A^\La.
\ee
Given the graded projection endomorphism 
\be
\vr=\op\sum_{k>0} \frac1k\ol\vr\circ h_k\circ h^n, \qquad
\ol\vr(\f)= \op\sum_{0\leq|\La|} (-1)^{\nm\La}\th^A\w
[d_\La(\dr^\La_A\rfloor\f)], \qquad \f\in \cS^{>0,n}_\infty[F;Y],
\ee
and the graded variational operator
$\dl=\vr\circ d$, the BGDA $\cS^*_\infty[F;Y]$
is split into the Grassmann-graded variational bicomplex
analogous to the above mentioned variational bicomplex of
$\cO^*_\infty Y$.
We restrict our consideration to its short variational subcomplex
\mar{g111}\beq
0\ar \Bbb R\ar \cS^0_\infty[F;Y]\ar^{d_H}\cS^{0,1}_\infty[F;Y]
\cdots \ar^{d_H} \cS^{0,n}_\infty[F;Y]\ar^\dl \bE_1, \quad
\bE_1=\vr(\cS^{1,n}_\infty[F;Y]), \label{g111}
\eeq
and subcomplex of one-contact graded forms
\mar{g112}\beq
 0\to \cS^{1,0}_\infty[F;Y]\ar^{d_H} \cS^{1,1}_\infty[F;Y]\cdots
\ar^{d_H}\cS^{1,n}_\infty[F;Y]\ar^\vr \bE_1\to 0. \label{g112}
\eeq
One can think of its elements
\mar{0709,'}\ben
&& L=\cL\om\in \cS^{0,n}_\infty[F;Y], \qquad \om=dx^1\w\cdots \w dx^n,
\label{0709}\\
&& \dl L=
\th^A\w
\cE_A\om=\op\sum_{0\leq|\La|}
 (-1)^{|\La|}\th^A\w d_\La (\dr^\La_A L)\om\in \bE_1 \label{0709'}
\een
as being a 
graded Lagrangian and its Euler--Lagrange operator, respectively.

Our goal now is cohomology of the subcomplexes (\ref{g111}) --
(\ref{g112}) of the BGDA $\cS^*_\infty[F;Y]$ and its de Rham complex
\mar{g110}\beq
0\to\Bbb R\ar \cS^0_\infty[F;Y]\ar^d \cS^1_\infty[F;Y]\cdots
\ar^d\cS^k_\infty[F;Y] \ar\cdots\,. \label{g110}
\eeq

\begin{theo} \label{v9} \mar{v9} There is an isomorphism
\mar{v10}\beq
H^*(\cS^*_\infty[F;Y])= H^*(Y) \label{v10}
\eeq
of cohomology $H^*(\cS^*_\infty[F;Y])$ of the
de Rham complex (\ref{g110}) to the de Rham cohomology  $H^*(Y)$ of $Y$. 
\end{theo} 

\begin{proof}
The complex (\ref{g110}) is the direct limit of the de Rham complexes of
the BGDAs $\cS^*[J^rY\op\times_X J^rF;J^rY]$, $r\in \Bbb N$.  Therefore,
the direct limit of cohomology groups
of these complexes is cohomology of the de Rham complex
(\ref{g110}). Cohomology of the de Rham complex of
$\cS^*[J^rY\op\times_X J^rF;J^rY]$ equals the de Rham
cohomology of
$J^rY$ \cite{bart,book05} and, consequently, that of $Y$, which is the
strong deformation retract of any $J^rY$. Hence,
the isomorphism (\ref{v10}) holds. 
\end{proof}

One can say something more. The isomorphism (\ref{v10}) is
induced by the cochain monomorphisms 
\be
\cO^*Y\to \cS^*[Y\op\times_X F;Y]\to \cS^*_\infty[F;Y].
\ee
Therefore,
any closed graded exterior form $\f\in
\cS^*_\infty[F;Y]$ is split into the sum $\f=d\si +\vf$ of an exact
graded exterior form  and a closed exterior form $\vf$ on $Y$.

Turn now to the complexes (\ref{g111}) --
(\ref{g112}). We have proved that,
in the case of an affine bundle $Y\to X$, cohomology of the short
variational complex (\ref{g111}) equals the de Rham cohomology of $X$,
while the complex (\ref{g112}) is exact \cite{cmp04}. Let
us generalize this result to the case of an arbitrary fiber bundle $Y\to
X$. 

\begin{lem} \label{v38'} \mar{v38'}
If $Y=\Bbb R^{n+m}\to \Bbb R^n$, the
complexes (\ref{g111}) -- (\ref{g112}) at all the terms, except $\Bbb R$,
are exact.
\end{lem}

\begin{proof}
This is the case of an affine bundle $Y$, and the above mentioned
exactness has been proved when the ring $\cO^0_\infty Y$ is restricted to
the subring $\cP^0_\infty Y$ of polynomial functions (see \cite{cmp04},
Lemmas 4.2 -- 4.3). The proof of these lemmas is straightforwardly
extended to
$\cO^0_\infty Y$ if the homotopy operator (4.5) in \cite{cmp04}, Lemma 4.2
is replaced with that (4.8) in \cite{cmp04}, Remark 4.1.
\end{proof}

\begin{theo} \label{v11} \mar{v11}
Cohomology of the complex (\ref{g111}) equals the de Rham cohomology
$H^*(Y)$ of $Y$. The complex (\ref{g112}) is exact. 
\end{theo}

\begin{proof}
The proof follows that of \cite{cmp04}, Theorem 2.1. 
We first prove Theorem \ref{v11} for the above mentioned BGDA
$\G(\gQ^*_\infty[F;Y])$. Similarly to
$\cS^*_\infty[F;Y]$, the sheaf
$\gQ^*_\infty[F;Y]$ and the BGDA $\G(\gQ^*_\infty[F;Y])$ are split into
the variational bicomplexes, and we consider their subcomplexes
\mar{v35-8}\ben
&& 0\ar \Bbb R\ar \gQ^0_\infty[F;Y]\ar^{d_H}\gQ^{0,1}_\infty[F;Y]
\cdots \ar^{d_H} \gQ^{0,n}_\infty[F;Y]\ar^\dl {\got E}_1, 
\label{v35}\\ 
&& 0\to \gQ^{1,0}_\infty[F;Y]\ar^{d_H}
\gQ^{1,1}_\infty[F;Y]\cdots
\ar^{d_H}\gQ^{1,n}_\infty[F;Y]\ar^\vr {\got E}_1\to 0, \label{v36}\\
&& 0\ar \Bbb R\ar \G(\gQ^0_\infty[F;Y])\ar^{d_H}\G(\gQ^{0,1}_\infty[F;Y])
\cdots \ar^{d_H} \G(\gQ^{0,n}_\infty[F;Y])\ar^\dl \G({\got E}_1),
\label{v37} \\
&&  0\to \G(\gQ^{1,0}_\infty[F;Y])\ar^{d_H}
\G(\gQ^{1,1}_\infty[F;Y])\cdots
\ar^{d_H}\G(\gQ^{1,n}_\infty[F;Y])\ar^\vr \G({\got E}_1)\to 0, \label{v38}
\een
where ${\got
E}_1 =\vr(\gQ^{1,n}_\infty[F;Y])$. By virtue of Lemma \ref{v38'}, the
complexes (\ref{v35}) -- (\ref{v36}) at all the terms, except $\Bbb R$, 
are exact. The terms $\gQ^{*,*}_\infty[F;Y]$ of
the complexes (\ref{v35}) -- (\ref{v36}) are sheaves of
$\G(\gQ^0_\infty)$-modules. Since $J^1Y$ admits a partition of
unity just by elements of $\G(\gQ^0_\infty)$, these sheaves are fine and,
consequently, acyclic. By virtue of the abstract de Rham theorem
(see \cite{cmp04}, Theorem 8.4, generalizing
\cite{hir}, Theorem 2.12.1), cohomology of the complex (\ref{v37}) equals
the cohomology of $J^\infty Y$ with coefficients in the constant sheaf
$\Bbb R$ and, consequently, the de Rham cohomology of $Y$, which is the
strong deformation retract of $J^\infty Y$. Similarly, the complex
(\ref{v38}) is proved to be exact. It remains to prove that cohomology of
the complexes (\ref{g111}) -- (\ref{g112}) equals that of the complexes
(\ref{v37}) -- (\ref{v38}). The proof of this fact straightforwardly
follows the proof of
\cite{cmp04}, Theorem 2.1, and it is a slight
modification of the proof of
\cite{cmp04}, Theorem 4.1, where graded exterior forms on the
infinite order jet manifold $J^\infty Y$ of an affine bundle are treated
as those on $X$. 
\end{proof}

\begin{prop} \label{cmp26} \mar{cmp26}
Every $d_H$-closed graded form $\f\in\cS^{0,m<n}_\infty[F;Y]$
falls into the sum
\mar{g214}\beq
\f=h_0\psi + d_H\xi, \qquad \xi\in \cS^{0,m-1}_\infty[F;Y], \label{g214}
\eeq
where $\psi$ is a closed $m$-form on $Y$. Every
$\dl$-closed graded Lagrangian $L\in \cS^{0,n}_\infty[F;Y]$ is the
sum
\mar{g215}\beq
\f=h_0\psi + d_H\xi, \qquad \xi\in \cS^{0,n-1}_\infty[F;Y], \label{g215}
\eeq
where $\psi$ is a closed $n$-form on $Y$.
\end{prop}

\begin{proof}
The complex (\ref{g111})
possesses the same cohomology as the similar part 
\mar{b317'}\beq
0\to\Bbb R\to \cO^0_\infty Y \ar^{d_H}\cO^{0,1}_\infty Y\cdots
\op\ar^{d_H} \cO^{0,n}_\infty Y \op\ar^\dl \bE_1 \label{b317'}
\eeq
of the variational
complex (\ref{b317}) of the GDA  $\cO^*_\infty Y$. The monomorphism
(\ref{v7}) and the body epimorphism (\ref{v7'}) yield the
corresponding cochain morphisms of the complexes (\ref{b317}) and
(\ref{b317'}). Therefore, cohomology of the complex (\ref{b317}) is the
image of cohomology of $\cO^*_\infty Y$.
\end{proof}

The global exactness of the complex (\ref{g112}) at the term
$\cS^{1,n}_\infty[F;Y]$ results in the following \cite{cmp04}.

\begin{prop} \label{g103} \mar{g103}
Given a graded Lagrangian $L=\cL\om$, there
is the decomposition
\mar{g99,'}\ben
&& dL=\dl L - d_H\Xi,
\qquad \Xi\in \cS^{1,n-1}_\infty[F;Y], \label{g99}\\
&& \Xi=\op\sum_{s=0}
\th^A_{\nu_s\ldots\nu_1}\w
F^{\la\nu_s\ldots\nu_1}_A\om_\la,\qquad
F_A^{\nu_k\ldots\nu_1}=
\dr_A^{\nu_k\ldots\nu_1}\cL-d_\la F_A^{\la\nu_k\ldots\nu_1}
+h_A^{\nu_k\ldots\nu_1},  \label{g99'}
\een
where local graded functions $h$ obey the relations
$h^\nu_a=0$,
$h_a^{(\nu_k\nu_{k-1})\ldots\nu_1}=0$.  Locally, one can
always choose $\Xi$ (\ref{g99'}) where all functions $h$ vanish.
\end{prop}

The decomposition (\ref{g99}  leads to the first variational formula for
graded Lagrangians as follows \cite{jmp05,cmp04}.
Let $\vt\in\gd \cS^0_\infty[F;Y]$ be a graded derivation of the $\Bbb
R$-ring $\cS^0_\infty[F;Y]$. The interior product $\vt\rfloor\f$ 
and the Lie
derivative $\bL_\vt\f$, $\f\in\cS^*_\infty[F;Y]$, are defined by
the formulae
\be
&& \vt\rfloor \f=\vt^\la\f_\la + (-1)^{[\f_A]}\vt^A\f_A, \qquad
\f\in \cS^1_\infty[F;Y],\\
&& \vt\rfloor(\f\w\si)=(\vt\rfloor \f)\w\si
+(-1)^{|\f|+[\f][\vt]}\f\w(\vt\rfloor\si), \qquad \f,\si\in
\cS^*_\infty[F;Y], \\
&& \bL_\vt\f=\vt\rfloor d\f+ d(\vt\rfloor\f), \qquad
\bL_\vt(\f\w\si)=\bL_\vt(\f)\w\si
+(-1)^{[\vt][\f]}\f\w\bL_\vt(\si).
\ee
A graded derivation $\vt$
 is said to be contact if the Lie
derivative $\bL_\vt$ preserves the ideal of contact graded forms
of the BGDA $\cS^*_\infty[F;Y]$. With respect to the local basis
$\{s^A\}$ for the BGDA $\cS^*_\infty[F;Y]$, any contact graded
derivation takes the form
\mar{g105}\beq
\vt=\vt_H+\vt_V=\vt^\la d_\la + (\vt^A\dr_A +\op\sum_{0<|\La|}
d_\La\vt^A\dr_A^\La), \label{g105}
\eeq
where the tuple of graded derivations $\{\dr_\la,\dr^\La_A\}$ is
defined as the dual of the tuple $\{dx^\la, ds^A_\La\}$ of
generating elements of the $\cS^0_\infty[F;Y]$-algebra
$\cS^*_\infty[F;Y]$, and $\vt^\la$, $\vt^A$ are local graded
functions \cite{cmp04}. One can justify that any vertical
contact graded derivation $\vt$ (\ref{g105}) satisfies the
relations
\mar{g232}\beq
\vt\rfloor d_H\f=-d_H(\vt\rfloor\f), \qquad
\bL_\vt(d_H\f)=d_H(\bL_\vt\f), \qquad \f\in\cS^*_\infty[F;Y].
\label{g232}
\eeq
Then it follows from the splitting (\ref{g99}) that the Lie derivative
$\bL_\vt L$ of a Lagrangian $L$ along a contact graded derivation
$\vt$ (\ref{g105}) fulfills the first variational formula
\mar{g107}\beq
\bL_\vt L= \vt_V\rfloor\dl L +d_H(h_0(\vt\rfloor \Xi_L)) + d_V
(\vt_H\rfloor\om)\cL, \label{g107}
\eeq
where $\Xi_L=\Xi+L$ is a Lepagean equivalent of $L$ given by the
coordinate expression (\ref{g99'}).

\end{document}